\documentclass{amsart}

\usepackage{amssymb}
\usepackage{amsmath}
\usepackage{amsthm}
\usepackage[mathscr]{eucal}

\newtheorem{thrm}{Theorem}[section]
\newtheorem{lemma}[thrm]{Lemma}
\newtheorem{prop}[thrm]{Proposition}
\newtheorem{cor}[thrm]{Corollary}

\theoremstyle{definition}
\newtheorem{defn}[thrm]{Definition}

\theoremstyle{remark}
\newtheorem{remark}[thrm]{Remark}

\numberwithin{equation}{section}

\begin{document}
\title{Analysis of the $\bar{\partial}$-Neumann problem along a straight edge}

\author{Dariush Ehsani}

\address{Department of Mathematics, Texas A\&M University, College
Station, Texas 77843-3368}
\email{ehsani@math.tamu.edu}

\subjclass[2000]{Primary 32W05; Secondary 35B65}

\begin{abstract} We show there exists an $L^p$ solution, for
$p\in(2,\infty)$, to the $\bar{\partial}$-Neumann problem on an
edge domain in $\mathbb{C}^2$ for $(0,1)$-forms, and we explicitly
compute the singularities, which are of complex logarithmic and
arctangent type, along the edge, of the solution.
\end{abstract}

\maketitle

\bibliographystyle{plain}
\section{Introduction}
The aim of this paper is to provide insight into the singular
behavior of solutions of the $\bar{\partial}$-Neumann problem on
domains which are not smooth.  We consider the edge domain,
$\Omega\subset\mathbb{C}^2$, defined by
\begin{equation*}
\{(z_1,z_2)\in\mathbb{C}^2:\Im z_1>\alpha\Im z_2>0\}
\end{equation*}
for some $0\leq\alpha<\infty$, and solutions to the
$\bar{\partial}$-Neumann problem on $\Omega$ for $(0,1)$-forms. A
solution to the $\bar{\partial}$-Neumann problem is an inverse to
the complex Laplacian,
$\bar{\partial}\bar{\partial}^{\ast}+\bar{\partial}^{\ast}\bar{\partial}$,
on $\Omega$.

The results obtained here are a generalization of results in
\cite{Eh}, which deals with the case of $\alpha=0$, in which the
singularities of the solution are explicitly computed.  Other
properties of the Neumann operator on non-smooth domains are
described in Ehsani \cite{Eh2}, Engli\v{s} \cite{Eng}, Henkin and
Iordan \cite{HI}, Henkin, Iordan, and Kohn \cite{HIK}, Michel and
Shaw \cite{MS1,MS2}, and Straube \cite{St}.

The domain of the edge considered here is an important model
domain in the study of the $\bar{\partial}$-Neumann problem on
non-smooth domains because, as in \cite{Eh} and \cite{Eh2}, we can
compute explicitly the singularities in the solution, however, on
the edge, the problem has the added complexity that the two
components $u_1$ and $u_2$ of the $(0,1)$-form solution
$u=u_1d\bar{z}_1+u_2d\bar{z}_2$ are coupled.  We resolve this
difficulty by examining the boundary conditions in detail along
the edge.  The domain is also important in that it depends on a
parameter, $\alpha$. Thus this domain should serve better as a
prototype for a wider class of non-smooth domains.

\section{Finding a solution}
\label{findsoln}
 We consider the $\bar{\partial}$-Neumann
problem on an edge, $\Omega$ in $\mathbb{C}^2$ described by
\begin{equation*}
\{(z_1,z_2)\in\mathbb{C}^2:\Im z_1>\alpha\Im z_2>0\}
\end{equation*}
for some $0\leq \alpha<\infty$.  The case of $\alpha=0$, in which
$\Omega$ is the cross product of two half-planes, was studied in
detail in \cite{Eh}.  For our data $(0,1)$-form, $f$, we make the
assumption $f\in\mathcal{S}_{(0,1)}(\overline{\Omega})$, the space
of $(0,1)$-forms whose coefficients are Schwartz functions. We use
the notation $z_j=x_j+iy_j$, for $j=1,2$. On the interior of
$\Omega$ the $\bar{\partial}$-Neumann problem becomes

\begin{equation*}
\triangle u_j=-2f_j \qquad j=1,2,
\end{equation*}
and the boundary conditions are
\begin{equation*}
u_2=0, \quad \frac{\partial u_1}{\partial \bar{z}_2}=0 \qquad
\mbox{on } y_2=0,
\end{equation*}
and
\begin{equation*}
u_1-\alpha u_2=0, \quad \frac{\partial u_1}{\partial
\bar{z}_2}-\frac{\partial u_2}{\partial \bar{z}_1}=0 \qquad
\mbox{on } y_1=\alpha y_2.
\end{equation*}
We make the change of coordinates
\begin{align*}
Y_2&=y_2\\
Y_1&=y_1-\alpha y_2,
\end{align*}
and we define the functions
\begin{align*}
u_{\alpha} &=u_1-\alpha u_2\\
f_{\alpha} &=f_1-\alpha f_2.
\end{align*}
In these new coordinates the interior equations become
\begin{equation}
\label{translap}
 \left(\frac{\partial^2}{\partial
x_1^2}+\frac{\partial^2}{\partial
x_2^2}+(1+\alpha^2)\frac{\partial^2}{\partial
Y_1^2}-2\alpha\frac{\partial^2}{\partial Y_1\partial
Y_2}+\frac{\partial^2}{\partial Y_2^2} \right)u_j=-2f_j
\end{equation}
for $j=\alpha,2$, on the product of two half-planes,
$\mathbb{H}\times\mathbb{H}=\{(x_1,x_2,Y_1,Y_2):Y_1,Y_2>0\}$, and
the boundary conditions become
\begin{align}
\label{dir1}
&u_{\alpha}=0 \quad \mbox{on } Y_1=0, \\
\label{dir2}
&u_2=0 \quad \mbox{on } Y_2=0,\\
\label{derY1} &-i\alpha\frac{\partial u_{\alpha}}{\partial
Y_1}=\frac{\partial u_2}{\partial x_1}-\alpha\frac{\partial
u_2}{\partial x_2}+i\left( (1+\alpha^2)\frac{\partial
u_2}{\partial Y_1}-\alpha\frac{\partial u_2}{\partial Y_2}\right)
 \quad \mbox{on } Y_1=0,\\
 \label{derY2}
&\frac{\partial u_{\alpha}}{\partial
x_2}+i\left(-\alpha\frac{\partial u_{\alpha}}{\partial
Y_1}+\frac{\partial u_{\alpha}}{\partial Y_2}\right)=-i
\alpha\frac{\partial u_2}{\partial Y_2}
 \quad \mbox{on } Y_2=0.
\end{align}

We apply the Fourier transform to (\ref{translap}) on the domain
$\mathbb{H}\times\mathbb{H}$.  We transform the equation for
$u_{\alpha}$.
\begin{multline}
\label{ftlapalpha}
-(\lambda_1^2+\lambda_2^2+(1+\alpha^2)\eta_1^2-2\alpha\eta_1\eta_2+\eta_2^2)
\hat{u}_{\alpha}\\
 -(1+\alpha^2)\left.\frac{\partial
\tilde{u}_{\alpha}}{\partial
Y_1}\right|_{Y_1=0}+i(2\alpha\eta_1-\eta_2)\tilde{\tilde{u}}_{\alpha}\Big|_{Y_2=0}-
\left.\frac{\partial \tilde{\tilde{u}}_{\alpha}}{\partial
Y_2}\right|_{Y_2=0}\\
=-2\hat{f}_{\alpha},
\end{multline}
where $\lambda_j$ is the transform variable corresponding to $x_j$
 and $\eta_j$ is the transform variable corresponding
to $Y_j$ for $j=1,2$, and $\tilde{u}_j$ denotes the partial
transform in all variables except $Y_1$ and $\tilde{\tilde{u}}_j$
denotes the partial transform of $u_j$ in all variables except
$Y_2$.

We use the superscript, $oj$, to denote an odd reflection with
respect to $Y_j$.
 Reflecting (\ref{ftlapalpha}) to be odd in
$\eta_1$, we have
\begin{multline}
\label{odd1}
-(\lambda_1^2+\lambda_2^2+(1+\alpha^2)\eta_1^2-2\alpha|\eta_1|\eta_2+\eta_2^2)
\hat{u}_{\alpha}^{o1}\\
+i(2\alpha|\eta_1|-\eta_2)\tilde{\tilde{u}}_{\alpha}^{o1}\Big|_{Y_2=0}-
\left.\frac{\partial \tilde{\tilde{u}}_{\alpha}^{o1}}{\partial
Y_2}\right|_{Y_2=0} =-2\hat{f}_{\alpha}^{o1}.
\end{multline}
We use (\ref{derY2}) to eliminate $\left.\frac{\partial
\tilde{\tilde{u}}_{\alpha}^{o1}}{\partial Y_2}\right|_{Y_2=0}$
from equation \ref{odd1}:
\begin{multline}
\label{alpha1}
-(\lambda_1^2+\lambda_2^2+(1+\alpha^2)\eta_1^2-2\alpha|\eta_1|\eta_2+\eta_2^2)
\hat{u}_{\alpha}^{o1}\\
+i(\alpha|\eta_1|-\eta_2-i\lambda_2)\tilde{\tilde{u}}_{\alpha}^{o1}\Big|_{Y_2=0}+\alpha
\left.\frac{\partial \tilde{\tilde{u}}_{2}^{o1}}{\partial
Y_2}\right|_{Y_2=0} =-2\hat{f}_{\alpha}^{o1}.
\end{multline}
We let $\zeta_1=\sqrt{\eta_1^2+\lambda^2}$ and set
$\eta_2=\alpha|\eta_1|-i\zeta_1$ in (\ref{alpha1}) in order to
eliminate $\tilde{\tilde{u}}_{\alpha}^{o1}\Big|_{Y_2=0}$. Finally
we solve for $\hat{u}_{\alpha}^{o1}$ in terms of
$\hat{f}_{\alpha}$ and $\left.\frac{\partial
\tilde{\tilde{u}}_2^{o1}}{\partial Y_2}\right|_{Y_2=0}$;
\begin{multline}
\label{alphawithder}
 \hat{u}_{\alpha}^{o1}=-i\alpha\frac{1}{\eta_2-\alpha|\eta_1|-i\zeta_1}
\frac{1}{\lambda_2-\zeta_1}\frac{\partial
\tilde{\tilde{u}}_{2}^{o1}}{\partial
Y_2}(\lambda_1,\lambda_2,\eta_1,0)\\
+\frac{2}{\lambda_1^2+\lambda_2^2+(1+\alpha^2)\eta_1^2-2\alpha|\eta_1|\eta_2+\eta_2^2}\times\\
\left(\hat{f}^{o1}_{\alpha}-\frac{\lambda_2+i(\alpha|\eta_1|-\eta_2)}{\lambda_2-\zeta_1}
\hat{f}_{\alpha}^{o1}(\lambda_1,\lambda_2,\eta_1,\alpha|\eta_1|-i\zeta_1)\right).
\end{multline}
Following an analogous procedure, we write
\begin{multline*}
 \hat{u}_2^{o2}=i\alpha\frac{1}{(1+\alpha^2)\eta_1-\alpha|\eta_2|+i\zeta_2}
 \frac{1}{\lambda_1-\alpha\lambda_2-\zeta_2}\frac{\partial
\tilde{u}_{\alpha}^{o2}}{\partial
Y_1}(\lambda_1,\lambda_2,0,\eta_2)\\+
\frac{2}{\lambda_1^2+\lambda_2^2+(1+\alpha^2)\eta_1^2-2\alpha\eta_1|\eta_2|+\eta_2^2}\times\\
\left(\hat{f}^{o2}_2-\frac{\lambda_1-\alpha\lambda_2+i\alpha|\eta_2|-i(1+\alpha^2)\eta_1}
{\lambda_1-\alpha\lambda_2-\zeta_2}
\hat{f}^{o2}_2(\lambda_1,\lambda_2,\frac{\alpha|\eta_2|-i\zeta_2}{1+\alpha^2},\eta_2)\right).
\end{multline*}

 An examination of the consistency
of the boundary conditions (\ref{derY1}) and (\ref{derY2}) along
the edge reveals
\begin{equation}
\left.\frac{\partial u_2}{\partial Y_2}\right|_{Y_1=Y_2=0}=
\left.\frac{\partial u_{\alpha}}{\partial Y_1}\right|_{Y_1=Y_2=0}.
\label{uniquechoice}
\end{equation}
Furthermore, relations (\ref{translap}), (\ref{dir1}), and
(\ref{dir2}) in (\ref{derY1}) when $Y_2=0$ allow us to determine
$\left.\frac{\partial^k}{\partial Y_1^k}\frac{\partial
u_2}{\partial Y_2}\right|_{Y_1=Y_2=0}$ ($k\geq 1$) in terms of
$\left.\frac{\partial^{k-1}}{\partial Y_1^{k-1}}\frac{\partial
u_2}{\partial Y_2}\right|_{Y_1=Y_2=0}$.  Similarly, relations
(\ref{translap}), (\ref{dir1}), and (\ref{dir2}) in (\ref{derY2})
when $Y_1=0$ give $\left.\frac{\partial^k}{\partial
Y_2^k}\frac{\partial u_{\alpha}}{\partial Y_1}\right|_{Y_1=Y_2=0}$
($k\geq 1$) in terms of $\left.\frac{\partial^{k-1}}{\partial
Y_2^{k-1}}\frac{\partial u_{\alpha}}{\partial
Y_2}\right|_{Y_1=Y_2=0}$. Thus, in the case $u\in
C^1_{(0,1)}(\overline{\Omega})$, $\left.\frac{\partial^k}{\partial
Y_1^k}\frac{\partial u_2}{\partial Y_2}\right|_{Y_1=Y_2=0}$ and
$\left.\frac{\partial^k}{\partial Y_2^k}\frac{\partial
u_{\alpha}}{\partial Y_1}\right|_{Y_1=Y_2=0}$ are finite $\forall
k \geq 0$.
\begin{remark}
\label{0term} By considering the decay of (\ref{alphawithder})
with respect to the Fourier variables, as they go to $\infty$,
from the condition $u\in C^1_{(0,1)}(\overline{\Omega})$ we can
conclude that
\begin{equation*}
\left.\frac{\partial u_2}{\partial Y_2}\right|_{Y_1=Y_2=0}=
\left.\frac{\partial u_{\alpha}}{\partial
Y_1}\right|_{Y_1=Y_2=0}=0
\end{equation*}
The last two terms of (\ref{alphawithder}) represent terms in
$C^1(\overline{\Omega})$, hence the decay of the first must be
sufficient enough to eliminate lower order terms (see \cite{Eh}
for details of this argument).
\end{remark}
In the case $u\in C^1_{(0,1)}(\overline{\Omega})$, the finite
Taylor coefficients, combined with the fact that
\begin{equation*}
\frac{\partial u_2}{\partial Y_2}(x_1,x_2,Y_1,0)\in
C^{\infty}(\mathbb{R}^2\times\mathbb{R}_+)
\end{equation*}
(see the arguments in Corollary \ref{reg} below), shows us
\begin{equation*}
\frac{\partial u_2}{\partial Y_2}(x_1,x_2,Y_1,0)\in
C^{\infty}(\mathbb{R}^2\times\overline{\mathbb{R}_+}).
\end{equation*}

From the symmetry of the domain in the $x$-variables and the fact
that $f\in\mathcal{S}_{(0,1)}(\overline{\Omega})$, we assume
$\left.\frac{\partial u_2}{\partial Y_2}\right|_{Y_2=0}$ is
Schwartz with respect to the $x$ variables, and, so that the
partial Fourier transform is determined (up to a
$C^{\infty}(\overline{\Omega})$ term) by the Taylor coefficients
at $Y_1=0$, we also assume $\left.\frac{\partial u_2}{\partial
Y_2}\right|_{Y_2=0}$ is Schwartz with respect to $Y_1$, and thus
that $\left.\frac{\partial u_2}{\partial Y_2}\right|_{Y_2=0}\in
\mathcal{S}(\mathbb{R}^2\times\overline{\mathbb{R}_+})$.

 We are therefore led to choose a $b_2(x_1,x_2,Y_1)\in
\mathcal{S}(\mathbb{R}^2\times\overline{\mathbb{R}_+})$ which
agrees to infinite order with $\frac{\partial u_2}{\partial
Y_2}(x_1,x_2,Y_1,0)$ at $Y_1=0$, possible by Borel's theorem.  We
also can choose in an analogous manner a
$b_{\alpha}(x_1,x_2,Y_2)\in
\mathcal{S}(\mathbb{R}^2\times\overline{\mathbb{R}_+})$ which
agrees to infinite order with $\frac{\partial u_{\alpha}}{\partial
Y_1}(x_1,x_2,0,Y_2)$ at $Y_2=0$.  The singular terms in the
solution we obtain are independent of the choice of $b_2$ and
$b_{\alpha}$, as the next lemma will show.
\begin{defn}
\label{equiv} We say $\hat{h}_1(\xi_1,\xi_2,\eta_1,\eta_2)$ is
equivalent to $\hat{h}_2(\xi_1,\xi_2,\eta_1,\eta_2)$, or
$\hat{h}_1 \sim \hat{h}_2$, if
$$h_1-h_2 \big|_{\Omega}\in C^{\infty}(\overline{\mathbb{H}\times\mathbb{H}}).$$
\end{defn}
\begin{lemma}
\label{db}
 Assume $u\in C^{1}_{(0,1)}(\overline{\Omega})$ and further
 that $\left.\frac{\partial u_2}{\partial Y_2}\right|_{Y_2=0}\in
\mathcal{S}(\mathbb{R}^2\times\overline{\mathbb{R}_+})$.
 Let $b_2\in
\mathcal{S}(\mathbb{R}^2\times\overline{\mathbb{R}_+})$ be chosen
as described above, then
\begin{multline}
\label{borel}
 \frac{1}{\eta_2-\alpha|\eta_1|-i\zeta_1}
\frac{1}{\lambda_2-\zeta_1} \frac{\partial
\tilde{\tilde{u}}_{2}^{o1}}{\partial
Y_2}(\lambda_1,\lambda_2,\eta_1,0)\sim\\
\frac{1}{\eta_2-\alpha|\eta_1|-i\zeta_1}
\frac{1}{\lambda_2-\zeta_1}
\tilde{\tilde{b}}_2^{o1}(\lambda_1,\lambda_2,\eta_1).
\end{multline}
  Also the relation, (\ref{borel}), is independent of the
choice of $b_2$.
\end{lemma}
In the proof of Lemma \ref{borel} we use the notation $\lesssim$
to mean $\leq c$ for $c>0$.
\begin{proof}
We first show
\begin{equation*}
\frac{1}{\eta_2-\alpha|\eta_1|-i\zeta_1}
\frac{1}{\lambda_2-\zeta_1}
\tilde{\tilde{b}}_2^{o1}(\lambda_1,\lambda_2,\eta_1)\in
L^p(\mathbb{R}^4)
\end{equation*}
for $p\in(1,2)$.  First integrating over $\eta_2$, we consider
\begin{multline*}
\int_{\mathbb{R}^4}\left|\frac{1}{\eta_2-\alpha|\eta_1|-i\zeta_1}
\frac{1}{\lambda_2-\zeta_1}
\tilde{\tilde{b}}_2^{o1}(\lambda_1,\lambda_2,\eta_1)
\right|^p d\mathbf{\lambda}d\mathbf{\eta} \lesssim\\
\int_{\mathbb{R}^3}\left|
\frac{\tilde{\tilde{b}}_2^{o1}(\lambda_1,\lambda_2,\eta_1)}
{\lambda_2-\zeta_1}\right|^p\frac{1}{\zeta_1^{p-1}}d\mathbf{\lambda}d\eta_1,
\end{multline*}
where $d\mathbf{\lambda}=d\lambda_1d\lambda_2$ and
$d\mathbf{\eta}=d\eta_1d\eta_2$.  Changing
$(\eta_1,\lambda_1,\lambda_2)$ to polar coordinates,
$(r,\phi,\theta)$, we then estimate
\begin{equation}
\label{firstint}
 \int_0^{2\pi}\int_0^{\pi}\int_0^{\infty}
\left|\frac{\tilde{\tilde{b}}_2^{o1}(\lambda_1,\lambda_2,\eta_1)}
{1-\cos\phi}\right|^p\frac{\sin\phi}{r^{2p-3}}drd\phi
d\theta.
\end{equation}
It is elementary to show, from the fact that
$b_2\in\mathcal{S}(\mathbb{R}^2\times\overline{\mathbb{R}}_+)$
and, from Remark \ref{0term}, which gives $b_2(Y_1=0)=0$, that
\begin{equation}
\label{derbnd}
 \left|\frac{\partial}{\partial
\eta_1}\tilde{\tilde{b}}_2^{o1}(\lambda_1,\lambda_2,\eta_1)\right|
\lesssim\frac{1}{1+r^3}.
\end{equation}

Therefore, with the Fundamental Theorem of Calculus, and
$\tilde{\tilde{b}}_2^{o1}(\lambda_1,\lambda_2,0)=0$,
(\ref{derbnd}) gives us the estimate
\begin{equation*}
\tilde{\tilde{b}}_2^{o1}(\lambda_1,\lambda_2,\eta_1)\lesssim\frac{1}{1+r^3}|\eta_1|,
\end{equation*}
 which, when used in
(\ref{firstint}), shows convergence of the integral.

 Now, if $v(x_1,x_2,Y_1)\in\mathcal{S}(\mathbb{R}^2\times\overline{\mathbb{R}}_+)$, and
$v$ vanishes to infinite order at $Y_1=0$ then after using a
partial Fourier inverse with respect to $\eta_2$ of
\begin{equation}
\label{vaninf} \frac{1}{\eta_2-\alpha|\eta_1|-i\zeta_1}
\frac{1}{\lambda_2-\zeta_1}
\tilde{\tilde{v}}^{o1}(\lambda_1,\lambda_2,\eta_1),
\end{equation}
we can use the decay of $v^{o1}(\lambda_1,\lambda_2,\eta_1)$,
faster than any power of $1/\eta_1$, to show (\ref{vaninf}) is
actually the transform of a function which, when restricted to
$\mathbb{H}\times\mathbb{H}$, is in
$C^{\infty}(\overline{\mathbb{H}\times\mathbb{H}})$. We denote by
$F.T._2$ the partial Fourier transform with respect to $Y_2$, and
$\Phi$ to be the Fourier inverse of (\ref{vaninf}).
\begin{align}
\nonumber
 |\zeta_1|^j\widehat{\frac{\partial^k}{\partial
Y_2^k}\Phi}&=
 |\zeta_1|^jF.T._2\left(\frac{\partial^k}{\partial
 Y_2^k}\tilde{\tilde{\Phi}}\right)\\
 \nonumber
 &=|\zeta_1|^j(i\alpha|\eta_1|-\zeta_1)^k\widehat{\Phi}\\
 \label{decay}
 &\lesssim
 \frac{1}{\eta_2-\alpha|\eta_1|-i\zeta_1}
\frac{1}{\lambda_2-\zeta_1}
|\zeta_1|^{j+k}\tilde{\tilde{v}}^{o1}(\lambda_1,\lambda_2,\eta_1).
\end{align}
Taking into account the decay of
$\tilde{\tilde{v}}^{o1}(\lambda_1,\lambda_2,\eta_1)$, we can show
(\ref{decay}) is in $L^p(\mathbb{R}^4)$ following the same proof
for $j=k=0$ above.

We then prove the lemma by setting
\begin{equation*}
v=\left.\frac{\partial u_2}{\partial Y_2}\right|_{Y_2=0}-b_2
\end{equation*}
above.
\end{proof}

As a corollary we have the
\begin{prop}
\label{solnlp}
 Let $u_{\alpha}$ and $u_2$ be defined on $\mathbb{H}\times\mathbb{H}$ in terms of
 their Fourier transforms as
\begin{multline}
\label{asoln}
\hat{u}_{\alpha}^{o1}=-i\alpha\frac{1}{\eta_2-\alpha|\eta_1|-i\zeta_1}
\frac{1}{\lambda_2-\zeta_1}\tilde{\tilde{b}}_{2}^{o1}(\lambda_1,\lambda_2,\eta_1)\\
+\frac{2}{\lambda_1^2+\lambda_2^2+(1+\alpha^2)\eta_1^2-2\alpha|\eta_1|\eta_2+\eta_2^2}\times\\
\left(\hat{f}^{o1}_{\alpha}-\frac{\lambda_2+i(\alpha|\eta_1|-\eta_2)}{\lambda_2-\zeta_1}
\hat{f}_{\alpha}^{o1}(\lambda_1,\lambda_2,\eta_1,\alpha|\eta_1|-i\zeta_1)\right)
\end{multline}
and
\begin{multline}
\label{2soln}
 \hat{u}_2^{o2}=i\alpha\frac{1}{(1+\alpha^2)\eta_1-\alpha|\eta_2|+i\zeta_2}
 \frac{1}{\lambda_1-\alpha\lambda_2-\zeta_2}
\tilde{b}_{\alpha}^{o2}(\lambda_1,\lambda_2,\eta_2)\\+
\frac{2}{\lambda_1^2+\lambda_2^2+(1+\alpha^2)\eta_1^2-2\alpha\eta_1|\eta_2|+\eta_2^2}\times\\
\left(\hat{f}^{o2}_2-\frac{\lambda_1-\alpha\lambda_2+i\alpha|\eta_2|-i(1+\alpha^2)\eta_1}
{\lambda_1-\alpha\lambda_2-\zeta_2}
\hat{f}^{o2}_2(\lambda_1,\lambda_2,\frac{\alpha|\eta_2|-i\zeta_2}{1+\alpha^2},\eta_2)\right).
\end{multline}
Then $u_{\alpha}$ and $u_2$ are in
$L^p(\mathbb{H}\times\mathbb{H})$ for $p\in(2,\infty)$.
\end{prop}
\begin{proof}
The first terms of the Fourier transforms, (\ref{asoln}) and
(\ref{2soln}), are in $L^p(\mathbb{R}^4)$ from the proof of Lemma
\ref{db}.  The proof that the last two terms in (\ref{asoln}) and
(\ref{2soln}) are in $L^p(\mathbb{R}^4)$ is the same as in the
case of $\alpha=0$ (see \cite{Eh}).  Then the Proposition follows
by the Hausdorff-Young theorem relating $L^p$ estimates of
functions in terms of $L^p$ estimates of their transforms.
\end{proof}
\begin{cor}
\label{reg} Let $u_{\alpha}$ and $u_2$ be defined as in
Proposition \ref{solnlp}.  Then $u_{\alpha}$ and $u_2$ are in
$C^{\infty}(\overline{V})$ for all neighborhoods
$V\subset\overline{\mathbb{H}\times\mathbb{H}}$ such that $V$ does
not intersect $\{Y_1=0\}\bigcap\{Y_2=0\}$.
\end{cor}
\begin{proof} We present the proof for $u_{\alpha}$.
Interior regularity follows from the strong ellipticity of the
Laplacian.

Also, general regularity at the boundary arguments for the
Dirichlet problem can be applied to the case in which $V$ is a
neighborhood such that $V\bigcap\partial
\left(\mathbb{H}\times\mathbb{H}\right)=V\bigcap\{Y_1=0\}\neq
\emptyset$ (see \cite{Fo}).

If $V$ is a neighborhood which intersects $Y_2=0$, then the
tangential derivatives commute with the $\bar{\partial}$-Neumann
problem in $V$, and thus as above, we can show $D_{\tau}^k
u_{\alpha}\in L^p(V)$ when $p>2$, for all tangential derivatives
$D_{\tau}^k$ of all orders $k$. Furthermore, since $u_2$ and
$\frac{\partial u_1}{\partial \bar{z}_2}-\frac{\partial
u_2}{\partial \bar{z}_1}$, in $\Omega$, satisfy Dirichlet
conditions along $y_2=0$, after a transformation, they belong to
$C^{\infty}(\overline{V})$, and hence we can also derive estimates
involving normal derivatives, and we conclude $D^ku_{\alpha}\in
L^p(V)$ when $p>2$ for all derivatives $D^k$ of all orders $k$.
Hence, a Sobolev embedding theorem applies to prove the corollary.
\end{proof}

With $u_{\alpha}$ and $u_2$ defined on
$\mathbb{H}\times\mathbb{H}$ as in Proposition \ref{solnlp}, and
with $u_1=u_{\alpha}+\alpha u_2$, we also denote by
$u_1(x_1,x_2,y_1,y_2)$ and $u_2(x_1,x_2,y_1,y_2)$ the
corresponding functions, defined on $\Omega$, under the
transformation $y_1=Y_1+\alpha Y_2$ and $y_2=Y_2$.

\begin{thrm}
\label{u1u2} With $u_{1}(x_1,x_2,y_1,y_2)$ and
$u_2(x_1,x_2,y_1,y_2)$ defined as above, the $(0,1)$-form,
$u=u_1d\bar{z}_1+u_2d\bar{z}_2$, is in
$C^1_{(0,1)}(\overline{\Omega})$ and $L^p_{(0,1)}(\Omega)$, for
$p\in(2,\infty)$, and solves the $\bar{\partial}$-Neumann problem
on the edge, $\Omega$, with data,
$f\in\mathcal{S}_{(0,1)}(\overline{\Omega})$.
 This solution is unique in the
sense that any other $C^1_{(0,1)}(\overline{\Omega})$ solution in
$L^p_{(0,1)}(\Omega)$ and whose boundary terms, $\frac{\partial
u_2}{\partial y_2}(y_2=0)$ and $\left(\frac{\partial u_1}{\partial
y_1}-\alpha\frac{\partial u_2}{\partial y_1}\right)(y_1=\alpha
y_2)$, are in
$\mathcal{S}(\mathbb{R}^2\times\overline{\mathbb{R}_+})$, differs
by a function in $C^{\infty}(\overline{\Omega})$.
\end{thrm}
\begin{proof} Remark \ref{0term} shows $u\in C_{(0,1)}^1(\overline{\Omega})$,
and Proposition \ref{solnlp} shows $u\in L^p_{(0,1)}(\Omega)$.
That $u$ solves the $\bar{\partial}$-Neumann problem follows by
our construction at the beginning of Section \ref{findsoln}.

The uniqueness part of the Proposition also follows from Remark
\ref{0term} which shows any solution in
$C^1_{(0,1)}(\overline{\Omega})$ and $L^p_{(0,1)}(\Omega)$ is
determined by choices of $b_2$ and $b_{\alpha}$, which are unique
modulo functions which vanish to infinite order at $Y_1=0$ and
$Y_2=0$, respectively, and thus, from Lemma \ref{db} the solutions
are unique modulo functions in $C^{\infty}(\overline{\Omega})$.
\end{proof}
\section{Singularities}
We shall examine the type of singularities which are present in
the solution described in Theorem \ref{u1u2}.  We shall proceed as
in \cite{Eh}, expanding $\hat{u}_{\alpha}^{o1}$ and
$\hat{u}_2^{o2}$ as asymptotic series for large $|\eta_1|$ and
$|\eta_2|$, in which higher order terms correspond to a class of
functions on $\mathbb{H}\times\mathbb{H}$ of greater
differentiability, continuous up to the boundary. We work with
$u_{\alpha}$, the analysis being similar for $u_2$. In what
follows, for $j=1,2$, let $\chi_{\eta_j} (\eta_j)$ be an even,
smooth function of $\eta_j$, with the property $\chi_{\eta_j}=1$
for $|\eta_j|<a$ and $\chi_{\eta_j}=0$ for $|\eta_j|>b$ for some
$b>a>0$.  Also, define $\chi_\eta (\eta_1,\eta_2)$ to be a smooth
function of $\eta_1$ and $\eta_2$, even in both variables, with
the property $\chi_\eta=1$ for $\eta_1^2+\eta_2^2<a$ and
$\chi_\eta=0$ for $\eta_1^2+\eta_2^2>b$ for some $b>a>0$.
\begin{lemma}  For $j=1,2$ let $\chi_{\eta_j}'=1-\chi_{\eta_j}$, and let $\chi_{\eta}'=1-\chi_{\eta}$.
With the equivalence relation defined in Definition \ref{equiv},
\begin{equation*}
\hat{u}_{\alpha}^{o1}\sim\chi_{\eta_1}'\chi_{\eta_2}'\hat{u}_{\alpha}^{o1}.
\end{equation*}
\end{lemma}
\begin{proof}[Sketch of proof]
  The
equivalence between $\hat{u}_{\alpha}^{o1}$ and
$\chi_{\eta}'\hat{u}_{\alpha}^{o1}$ is obvious, and that between
$\chi_{\eta}'\hat{u}_{\alpha}^{o1}$ and
$\chi_{\eta_1}'\chi_{\eta_2}'\hat{u}_{\alpha}^{o1}$ may be shown
by evaluating decay properties of
 $(\chi_{\eta}'-\chi_{\eta_1}'\chi_{\eta_2}')\hat{u}_{\alpha}^{o1}$
 in Fourier transform space.
\begin{equation}
\label{chi}
 \chi_{\eta}'-\chi_{\eta_1}'\chi_{\eta_2}'=\left(\chi_{\eta_1}-\chi_{\eta}\right)
 +\chi_{\eta_2}\left(1-\chi_{\eta_1}\right).
\end{equation}
When the first term on the right hand side of (\ref{chi}) is
multiplied by each term of $\hat{u}_{\alpha}^{o1}$, as expressed
by (\ref{asoln}), after taking a partial Fourier inverse with
respect to $\eta_2$, we can use the relation between taking
derivatives with respect to $Y_2$ and multiplying by $\zeta_1$ as
in (\ref{decay}) to show differentiability in all variables given
the decay with respect to the Fourier variables $\lambda_1$,
$\lambda_2$ and $\eta_1$.

When the second term on the right hand side of (\ref{chi}) is
multiplied by each term of $\hat{u}_{\alpha}^{o1}$, as expressed
by (\ref{asoln}), we can again use the relation between taking
derivatives with respect to $Y_2$ and multiplying by $\zeta_1$,
this time using decay with respect to $\eta_2$ to derive decay
with respect to $\eta_1$ to finish the proof of the lemma.
\end{proof}

To obtain our asymptotic expansion of
$\chi_{\eta_1}'\chi_{\eta_2}'\hat{u}_{\alpha}^{o1}$, we expand
$\zeta_1$ for large $|\eta_1|$,
\begin{equation*}
\frac{1}{\lambda_1^2+\lambda_2^2+(1+\alpha^2)\eta_1^2-2\alpha|\eta_1|\eta_2+\eta_2^2}
=\frac{1}{\lambda_1^2+\lambda_2^2+(\eta_2-\alpha|\eta_1|)^2+\eta_1^2}
\end{equation*}
as a geometric series in $(\eta_2-\alpha|\eta_1|)^2+\eta_1^2$, and
we integrate by parts all Fourier integrals involving $f_{\alpha}$
or $b_2$, leaving as remainders those terms which decay faster
than either of
\begin{equation*}
\frac{1}{|\eta_1|^{2n+3}}
\frac{1}{\lambda_1^2+\lambda_2^2+(\eta_2-\alpha|\eta_1|)^2+\eta_1^2}
\end{equation*}
or
\begin{equation*}
\frac{1}{|\eta_2|^{2n+3}}
\frac{1}{\lambda_1^2+\lambda_2^2+(\eta_2-\alpha|\eta_1|)^2+\eta_1^2}
\end{equation*}
for large $|\eta_1|$ or $|\eta_2|$. Again relating $Y_2$
derivatives of the partial Fourier inverse of
$\frac{1}{\lambda_1^2+\lambda_2^2+(\eta_2-\alpha|\eta_1|)^2+\eta_1^2}$
with multiplication by $\zeta_1$, we can show decay in $|\eta_1|$
gives differentiability with respect to $Y_2$ and vice-versa, and
thus all remainder terms are Fourier transforms of functions in
$C^{n}(\overline{\mathbb{H}\times\mathbb{H}})$.

Our asymptotic expansion, for large $|\eta_1|$, $|\eta_2|$, is
thus a sum of terms of the form
\begin{equation}
\label{asexpn}
\chi_{\eta_1}'\chi_{\eta_2}'c_{jklm}(\lambda_1,\lambda_2)
\frac{1}{\eta_2^j}\frac{1}{|\eta_1|^k}\frac{1}{\eta_1^{2l-1}}
\frac{1}{((\eta_2-\alpha|\eta_1|)^2+\eta_1^2)^m}
\end{equation}
for $j\geq -1$, $k,m=0,1$, and $l,n\geq 1$, where
$c_{jklm}(\lambda_1,\lambda_2)$ are in
$\mathcal{S}(\mathbb{R}^2)$.

 We start with the terms
\begin{equation*}
\frac{\chi_{\eta_1}'\chi_{\eta_2}'}{\left((\eta_2-\alpha|\eta_1|)^2+\eta_1^2\right)^{j+1}}.
\end{equation*}
For $\eta_1\neq 0$ and $Y_2>0$
\begin{multline*}
\int_{-\infty}^{\infty}\frac{1}{\left((\eta_2-\alpha|\eta_1|)^2+\eta_1^2\right)^{j+1}}
e^{i\eta_2Y_2} d\eta_2=\\ \frac{2\pi i}{j!}\sum_{k=0}^j {j\choose
k}
(-1)^{j-k}\frac{(2j-k)!}{j!}\frac{(iY_2)^k}{(2i|\eta_1|)^{2j-k+1}}e^{|\eta_1|(-1+i\alpha)Y_2},
\end{multline*}
which is a linear combination of terms of the form
\begin{equation*}
\int_0^{Y_2}\cdots\int_0^{t_2}
\int_{-\infty}^{\infty}\frac{1}{\left(\eta_2^2+\eta_1^2\right)^{l}}
e^{i\eta_2(1-i\alpha)Y_2} d\eta_2dt_1\cdots dt_{j+1-l}.
\end{equation*}
Such terms (excluding the constants of integration, whose inverses
are singular along all of $Y_2>0$) were studied in Lemma 3.7 of
\cite{Eh}.  And from the same Lemma 3.7, which gives
$\frac{\chi_{\eta_1}'\chi_{\eta_2}'}{\eta_1^2+\eta_2^2}$ locally
near $Y_1=Y_2=0$, we immediately have
\begin{lemma}
The inverse Fourier transform of
\begin{equation*}
\frac{\chi_{\eta_1}'\chi_{\eta_2}'}
{\left((\eta_2-\alpha|\eta_1|)^2+\eta_1^2\right)^{j+1}},
\end{equation*}
near $Y_1=Y_2=0$, has the form
\begin{equation}
\label{logform}
 p(Y_1,Y_2)\log(Y_1^2+(1-i\alpha)^2Y_2^2),
\end{equation}
where $p$ is a homogeneous polynomial of degree $2j$ in $Y_1$ and
$Y_2$, modulo functions which are in
$C^{\infty}(\overline{\mathbb{H}\times\mathbb{H}})$ or are
singular along all of $Y_2>0$.
\end{lemma}
With a slight abuse of notation we shall use the equivalence
relation in Definition \ref{equiv} to apply to functions defined
on $\mathbb{H}\times\mathbb{H}$.

We now define functions $\Phi_l$ on $Y_2\geq 0$ which have the
form of (\ref{logform}) such that
\begin{equation*}
\chi_{\eta_1}'\chi_{\eta_2}'\widehat{\chi\Phi_l}\sim
\frac{\chi_{\eta_1}'\chi_{\eta_2}'}
{{\left((\eta_2-\alpha|\eta_1|)^2+\eta_1^2\right)^{l}}}.
\end{equation*}
Then with $\Phi_l$ defined for $l\geq 1$, we define $(\Phi_l)_0=
\Phi_l$ for $Y_2\geq 0$, and, for $j \geq 1$, $(\Phi_l)_j$ to be
the unique solution of the form
\begin{equation*}
 p_1\log (Y_1^2+(1-i\alpha)^2Y_2^2) + p_2 +p_3\arctan
\left(\frac{Y_1}{(1-i\alpha)Y_2}\right)
\end{equation*}
 on the half-plane
$\{(Y_1,Y_2):Y_2\geq 0\}$, where $p_1$, $p_2$, and $p_3$ are
polynomials in $Y_1$ and $Y_2$ such that $p_2(0,Y_2)=0$, to the
equation
\begin{equation*}
\frac{\partial(\Phi_l)_j}{\partial Y_1}=(\Phi_l)_{j-1}.
\end{equation*}
Also, define for $k\geq 1$, on $Y_2\geq 0$ and restricting to
$Y_1\geq 0$,
\begin{equation*}
(\Phi_l)_{jk}=\int_0^{Y_2}\cdots\int_0^{t_2}\int_0^{t_1}
(\Phi_l)_j(Y_1,t)dt dt_1\cdots dt_{k-1}.
\end{equation*}
Then integration by parts in the Fourier transform integral shows
\begin{equation*}
\chi_{\eta_1}'\chi_{\eta_2}'\widehat{\chi(\Phi_l)_{jk}^{o1}}\sim
\frac{\chi_{\eta_1}'\chi_{\eta_2}'}{|\eta_1|^m\eta_1^{2n+1}\eta_2^k}
\frac{1}{{\left((\eta_2-\alpha|\eta_1|)^2+\eta_1^2\right)^{l}}},
\end{equation*}
where $2n+1+m=j$.

We are now ready to prove the
\begin{thrm}
\label{mnthrm}
 Let $f\in\mathcal{S}_{(0,1)}(\overline{\Omega})$,
and $u=u_1d\bar{z}_1+u_2d\bar{z}_2$ be the $(0,1)$-form which
solves the $\bar{\partial}$-Neumann problem on $\Omega$ with data
$f$. Then, in $\Omega$, near $y_1=y_2=0$, $u_j$ can be written as
\begin{multline}
\label{formofu}
u_j=\\
\alpha_{j1}\log((y_1-\alpha
y_2)^2+(1-i\alpha)^2y_2^2)+\alpha_{j2}\log((y_1-\alpha
y_2)^2(1-i\alpha)^2+y_2^2)\\
+\beta_{j1}\arctan\left(\frac{y_1-\alpha
y_2}{(1-i\alpha)y_2}\right)+\beta_{j2}\arctan\left(
\frac{(1-i\alpha)(y_1-\alpha y_2)}{y_2}\right)+\gamma_j,
\end{multline}
where $\alpha_{jk}$, $\beta_{jk}$ and $\gamma_j$ are smooth for
$j,k=1,2$.
\end{thrm}
\begin{proof}
We may use the functions $(\Phi_l)_{jk}$ constructed above, which
have the form
\begin{equation*}
(\Phi_l)_{jk}= p_1\log (Y_1^2+(1-i\alpha)^2Y_2^2) + p_2
+p_3\arctan \left(\frac{Y_1}{(1-i\alpha)Y_2}\right) +p_4\log
|Y_1|,
\end{equation*}
where the $p_m$ are homogeneous polynomials of degree $(2l-2)+j+k$
in $Y_1$ and $Y_2$ for $m=1,2,3,4$,
 to see the structure of the terms of the
form
\begin{equation*}
\frac{\chi_{\eta_1}'\chi_{\eta_2}'}{|\eta_1|^m\eta_1^{2n+1}\eta_2^k}
\frac{1}{{\left((\eta_2-\alpha|\eta_1|)^2+\eta_1^2\right)^{l}}}
\end{equation*}
arising in the asymptotic expansion for $u_{\alpha}$.  For the
other terms, of the form,
\begin{equation*}
\chi_{\eta_1}'\chi_{\eta_2}'\frac{\eta_2}{|\eta_1|^m\eta_1^{2n+1}}
\frac{1}{{\left((\eta_2-\alpha|\eta_1|)^2+\eta_1^2\right)^{l}}}
\end{equation*}
we use the property
\begin{equation*}
\frac{\partial \Phi_l}{\partial
Y_2}=\frac{(1-i\alpha)^2}{2(j-1)}y_2\Phi_{l-1}.
\end{equation*}
Since Corollary \ref{reg} allows us to conclude any singular terms
along all of $Y_1=0$ or $Y_2=0$ must vanish, we can take a finite
number of terms of the form (\ref{asexpn}) in the asymptotic
expansion and pair each with an appropriate function constructed
with the $(\Phi_l)_{jk}$, ignoring singular terms such as
$\log|y_1|$ to show $\forall n \in \mathbb{N}$, $\exists$
polynomials, $A_n$, $B_n$, and $C_n$, of degree $n$ in $Y_1$ and
$Y_2$, and whose coefficients are Schwartz functions of
$\lambda_1$ and $\lambda_2$, and $D_n$, the partial transform in
the $x$ variables of a function which belongs to
$C^{n}(\overline{\mathbb{H}\times\mathbb{H}})$, such that near
$Y_1,Y_2=0$
\begin{multline*}
  F.T._x \left( u_\alpha^{o1}
\right)(\lambda_1,\lambda_2,Y_1,Y_2) =\\
A_n \log (Y_1^2+(1-i\alpha)^2Y_2^2) + B_n +C_n \arctan \left(
\frac{Y_1}{(1-i\alpha)Y_2} \right) + D_n,
\end{multline*}
where $F.T._x$ stands for the partial Fourier transform in the $x$
variables.

Lastly, using Borel's theorem, inverting with respect to
$\lambda_1$ and $\lambda_2$, and transforming back to the
variables, $y_1$ and $y_2$, we can show $u_{\alpha}$ is of the
form (\ref{formofu}).  Then combining with an analogous argument
applied to $u_2$, we conclude the theorem.
\end{proof}

We end with the note that Theorem \ref{mnthrm} is non-trivial;
there exists an $f\in\mathcal{S}_{(0,1)}(\overline{\Omega})$, for
instance an $f\in \mathcal{S}_{(0,1)}(\overline{\Omega})$ which is
equivalently equal to 1 in a neighborhood of the edge, such that
one of the $\alpha_{ij}$ or $\beta_{ij}$ is not equivalently 0.

\end{document}